%
%
%
%
%
%
%
%
%
%
\scrollmode
\magnification=\magstep1
\parskip=\smallskipamount

\def\demo#1:{\par\medskip\noindent\it{#1}. \rm}
\def\ni{\noindent}               
\def\ll{\leftline}
\def\cl{\centerline}

%
%
\outer\def\beginsection#1\par{\bigskip
  \message{#1}\leftline{\bf\&#1}
  \nobreak\smallskip\vskip-\parskip\noindent}

%
%
\outer\def\proclaim#1:#2\par{\medbreak\vskip-\parskip
    \noindent{\bf#1.\enspace}{\sl#2}
  \ifdim\lastskip<\medskipamount \removelastskip\penalty55\medskip\fi}

\def\endpr{\hfill $\spadesuit$ \medskip}

%
%
%
%


%
%
%
%
\def\B{{\bf B}}
\def\C{{\bf C}}

\def\N{{\bf N}}


%
%
%
%

%
%
%

\def\d{\delta}
\def\e{\epsilon}
\def\z{\zeta}

\def\x{\xi}


%
%
%
%
\def\bs{\backslash}              

%
%
\def\holo{holomorphic}                   

\def\pc{polynomially convex}             

\def\dist{{\rm dist}}

\def\Aut{{\rm Aut}}                         

\def\begin{\ll{} \vskip 10mm \nopagenumbers}  

\begin
\cl{\bf AN INTERPOLATION THEOREM FOR}
\medskip
\cl{\bf HOLOMORPHIC AUTOMORPHISMS OF $\C^{\bf n}$}
\bigskip\medskip
\cl{\bf Gregery T. Buzzard\footnote{*}
{\rm Research at MSRI is supported in part by NSF grant DMS-9022140}
 and Franc Forstneric\footnote{**}
{\rm Supported in part by an NSF grant and by a grant from the Ministry
of Science of the Republic of Slovenia}}
\bigskip\rm
\ll{Department of Mathematics, Indiana University, Bloomington, IN 47405, USA}
\ll{Department of Mathematics, University of Wisconsin, Madison, WI 53706, USA}
\medskip\bigskip

\font\small=cmr9

\midinsert
\narrower\narrower
\cl{\bf Abstract} 
\baselineskip=10pt
{\small We construct automorphisms of $\C^n$ which map
certain discrete sequences one onto another with prescribed
finite jet at each point, thus solving a general Mittag-Leffler
interpolation problem for automorphisms.  Under certain circumstances,
this can be done while also approximating a given automorphism on a
compact set.}
\endinsert

\baselineskip=\normalbaselineskip
\beginsection 0. Introduction

Let $\C^n$ be complex Euclidean space of dimension $n$ with
coordinates $z=(z_1, \ldots, z_n)$. We shall always assume $n>1$.
We denote by $\Aut\C^n$ the group of all holomorphic automorphisms
of $\C^n$, and by $\Aut_1\C^n$ the group of automorphisms $F\in\Aut\C^n$
with complex Jacobian one: $JF\equiv 1$. For some recent
developments regarding these groups see the survey [2].
By $\B$ we denote the open unit ball in $\C^n$.

In this paper we construct automorphisms of $\C^n$ which map
certain discrete sequences one onto another with prescribed
finite jet at each point, thus solving a general Mittag-Leffler
interpolation problem for automorphisms.

In approaching such an interpolation problem one must consider
the limitations of automorphisms in mapping one sequence
to another.  Rosay and Rudin [5] showed that, in general, one cannot map
an infinite discrete set in $\C^n$ onto another discrete set
by an automorphism of $\C^n$. In fact, the infinite discrete sets
in $\C^n$ form uncountably many different equivalence classes under
this relation [5, sect. 5]. They introduced the following notion:

\proclaim Definition: A discrete infinite sequence
$\{a_j\} \subset\C^n$ (without repetition) is {\it tame} if
there exists a \holo\ automorphism $F\in\Aut\C^n$ such that
$$ F(a_j)=je_1=(j,0,\ldots,0),\qquad j=1,2,3,\ldots  $$
The sequence $\{a_j\}$ is {\it very tame\/} if the above holds
for some $F\in\Aut_1\C^n$.

Rosay and Rudin [5] showed that points in a tame sequence are
permutable by automorphisms of $\C^n$.
Hence one can speak about the tameness of an infinite discrete
set. On the other hand, they showed that there exist discrete sets
$E\subset \C^n$ such that the only nondegenerate \holo\ map
$F\colon \C^n\to\C^n$ satisfying $F(\C^n\bs E)\subset \C^n\bs E$
is the identity map [5, theorem 5.1].

In view of this, it is reasonable to consider tame sequences
in order to get positive results on jet interpolation.
The main result of this paper is the following.

%
%
%
\proclaim 0.1 Theorem:
Assume that $n>1$, that $\{a_j\}$ and $\{b_j\}$ are tame
sequences in $\C^n$ (without repetitions), and that $P_j\colon\C^n \to \C^n$
is a \holo\ polynomial map of degree at most $m_j \ge 1$ such that
$P_j(0) = 0$ and $JP_j(0) \ne 0$.  Then there exists $F \in
\Aut\C^n$ such that for every $j = 1,2,3,\ldots$ we have $F(a_j)=b_j$ and
$$F(z)=b_j+P_j(z-a_j)+O(|z-a_j|^{m_j+1}), \quad  z\to a_j. \eqno(0.1)$$
If in addition the polynomials $P_j$ satisfy
$$ JP_j(z)=1+O(|z|^{m_j}),\ \ z\to 0,\quad j=1,2,3,\ldots, \eqno(0.2)$$
and the sequences $\{a_j\}$ and $\{b_j\}$ are very tame, we may choose
$F \in \Aut_1\C^n$.   

A version of this theorem in which the degrees of the polynomials $P_j$
are uniformly bounded can be found in [3, theorem 3.1].  In that
result, with $a_j=b_j=je_1$ for all $j$, the automorphism $F$ is a
{\it finite} composition of shears and generalized shears.
See the survey [2] for further information on shears and
generalized shears.

Under certain circumstances we can obtain the above result while
at the same time approximating a given automorphism on a \pc\ set.

\footline={\hss\tenrm\folio\hss}

\proclaim 0.2 Theorem: In addition to the hypotheses of
theorem 0.1, assume that $\Phi\in\Aut\C^n$ is an automorphism
and $K\subset \C^n\bs \{a_j\}_{j=1}^\infty$ is a compact, \pc\ set
such that $\Phi(K) \subset \C^n\bs \{b_j\}_{j=1}^\infty$.
Then for each $\e>0$ there exists an $F\in\Aut\C^n$
satisfying theorem 0.1 and
$$ |F(z)-\Phi(z)|<\e,\qquad z\in K.                 \eqno(0.3) $$
If the volume preserving assumptions in theorem 0.1 hold
and if $\Phi\in\Aut_1\C^n$, we may choose $F\in\Aut_1\C^n$.

\medskip
The paper is organized as follows.
In section\ 1 we recall a result concerning holomorphic
automorphisms of $\C^n$ which have prescribed jets of finite order at
a finite set of points, and which are close to the identity on a given
\pc\ set, as well as a result about convergence of compositions of
automorphisms.  Both of these results are proved in [3].

In section\ 2 we provide the proof of theorem 0.1.  The outline of the
proof is as follows.  Since $\{a_j\}$ and $\{b_j\}$ are tame, we may
first assume that they both lie in the $z_1$-axis.  Next, if $\{c_j\}$
is a tame set in the $z_1$-axis, we can find $\Psi \in \Aut_1 \C^n$ such
that at each point $a_j$, $\Psi$ agrees with the translation $z \mapsto
c_j + (z-a_j)$ to order $m_j$.  In fact, we can find such a $\Psi$
which is a composition of three shears of the form $z \mapsto z +
f(\lambda(z)) v$, where $f$ is entire on $\C$, $v \in \C^n\bs\{0\}$, and
$\lambda$ is a $\C$-linear form with $\lambda(v) = 0$.  

Hence it suffices to solve
the following simpler problem:  Given polynomials $P_j$ as above, 
construct a tame sequence $\{c_j\}$ contained in the
$z_1$-axis and $F \in \Aut\C^n$ such that $F$ fixes each $c_j$ and has
prescribed jet $P_j$ at $c_j$ (with appropriate modifications in the
volume-preserving case). 

To construct $F$ and $\{c_j\}$, we use an inductive procedure.  We
take $c_1=3e_1$ and use the results of section 1 to find an
automorphism $F_1$ with prescribed jet $P_1$ at $c_1$ and which is
near the identity on $K_1 = 2 \overline{\B}$, where $\B$ is the unit
ball in $\C^n$.  We then choose a point $c_2$ on the $z_1$-axis so
that $c_2$ and $F_1(c_2)$ lie outside the convex hull $K_2$ of
$2\overline{\B} \cup K_1 \cup F_1(2\overline{\B}) \cup \{c_1\}$.  Again
the results of section 1 allow us to find an automorphism $\Psi_2$
which is near the identity on $K_2$ and which agrees with the identity
to order $m_1$ at $c_1$, which maps $F_1(c_2)$ to $c_2$, and which has
an appropriate jet at $F_1(c_2)$.  Then $F_2 = \Psi_2 \circ F_1$.
We then repeat this procedure, choosing $c_{k+1}$ so that $c_{k+1}$
and $F_k(c_{k+1})$ lie outside the convex hull $K_{k+1}$ of 
$(k+1)\overline{\B} \cup K_k \cup F_k((k+1)\overline{\B}) \cup \{c_k\}$,
then choosing $\Psi_k$ to be near the identity on $K_{k+1}$, 
equal to the identity to order $m_j$ at $c_j$ for $j \le k$, mapping
$F(c_{k+1})$ to $c_{k+1}$ and having appropriate jet at $F_k(c_{k+1})$.  

In this way we define a sequence of automorphisms $F_1, F_2, \ldots$
with $F_{k+1} = \Psi_{k+1} \circ F_k$.  Then the second result of section 1
implies that this sequence converges to an automorphism of $\C^n$ with
the desired properties. 

In section 3 we prove theorem 0.2.  We use the result of theorem 0.1,
but must take extra care to insure that we can in addition
approximate the given map.  The main addition here is lemma 3.1:
given a compact, \pc\ set $K\subset \C^n$, an automorphism
$H$ of $\C^n$, and a discrete sequence $\{a_j\}$ contained in
the $z_1$-axis outside $K$, there is an $R>0$ and an automorphism
$\Psi\in\Aut\C^n$ that approximates $H$ on $K$ and satisfies
$\Psi(a_j)=(R+j)e_1$ for each $j$.

%
%
%
\beginsection 1.  Interpolation on a finite set and convergence of
automorphisms

The results of this section are taken directly from [3].
We do not reproduce the proofs here. The first result concerns
jet interpolation by automorphisms on a finite set.


\medskip\ni\bf 1.1 Proposition. \sl Let $n>1$. Assume that
\item{(a)} $K\subset\C^n$ is a compact, \pc\ set,
\item{(b)} $\{a_j\}_{j=1}^k \subset K$ is a finite subset of $K$,
\item{(c)} $p$ and $q$ are arbitrary points in $\C^n\bs K$ (not
necessarily distinct),
\item{(d)} $N$ is a nonnegative integer, and
\item{(e)} $P\colon \C^n\to \C^n$ is a holomorphic polynomial map
of degree at most $m\ge 1$, satisfying $P(0)=0$ and $JP(0)\ne 0$.

\ni Then for each $\e>0$ there exists an automorphism
$F\in \Aut\C^n$ satisfying
\item{(i)}   $F(z)=q+P(z-p)+O(|z-p|^{m+1})$ as $z\to p$,
\item{(ii)}  $F(z)=z+O(|z-a_j|^N)$ as $z\to a_j$ for each
$j=1,2,\ldots,k$, and
\item{(iii)} $|F(z)-z|+|F^{-1}(z)-z|<\e$ for each $z\in K$.

\ni If in addition the polynomial map $P$ satisfies
\item{(e')} $P(0)=0$ and $JP(z)=1+O(|z|^m)$ as $z\to 0$,

\ni then there exists a polynomial automorphism $F$ with
$JF=1$ satisfying (i)--(iii).
\medskip\rm

The next result concerns convergence of a sequence of automorphisms.
It gives sufficient conditions for such a sequence to converge to a
biholomorphic map in some domain and describes the domain of
convergence and the image of the limit map.

\proclaim 1.2 Proposition:
Let $D$ be a connected open set in $\C^n$ which is exhausted
by compact sets
$K_0\subset K_1\subset K_2\subset\cdots \subset \cup_{j=0}^\infty K_j =D$
such that $K_{j-1} \subset {\rm Int} K_j$ for each $j\in\N$.
Choose numbers $\e_j$ ($j=1,2,3,\ldots$) such that
$$ 0< \e_j < \dist(K_{j-1},\C^n\bs K_j)\ \ (j\in\N), \qquad
      \sum_{j=1}^\infty \e_j <\infty.      \eqno(1.1) $$
Suppose that for each $j=1,2,3,\ldots$, $\Psi_j$ is a holomorphic
automorphism of $\C^n$ satisfying
$$ |\Psi_j(z)-z| < \e_j, \quad z\in K_j.           \eqno(1.2) $$
Set $\Phi_m=\Psi_m\circ\Psi_{m-1}\circ\cdots \circ \Psi_1$.
Then there is an open set $\Omega\subset\C^n$ such that
$\lim_{m\to\infty} \Phi_m =\Phi$ exists on $\Omega$
(uniformly on compacts), and $\Phi$ is a biholomorphic map of
$\Omega$ onto $D$.  In fact, $\Omega = \cup_{m=1}^\infty
\Phi_m^{-1}(K_m)$.  

{}Finally, we include some elementary observations about power series.
If $F\colon U\subset\C^n \to\C^n$ is a \holo\ map, $a\in U$,
and $m\ge 1$ an integer, we write
$$ F(z)=F(a)+ F_{m,a}(z-a) + O(|z-a|^{m+1}),\quad z\to a.  $$
Thus $F_{m,a}$ is just the Taylor polynomial of $F$ of order $m$
at $a$ without the constant term. 

The following lemma is evident by composing the power series.

\proclaim 1.3 Lemma: If $F$ and $G$ are nondegenerate \holo\ maps
on certain open subsets of $\C^n$, with values in $\C^n$, then
for each integer $m\ge 1$ we have
$$ (G\circ F)_{m,a}(z) =
   G_{m,F(a)} \circ F_{m,a}(z)+ O(|z|^{m+1}),\quad z\to 0.    $$
at each point $a$ where the composition $G\circ F$ is defined.

This generalizes to composition of several maps. We also have
$$ (F^{-1})_{m,F(a)} \circ F_{m,a}(z)=z+O(|z|^{m+1}), \quad z\to 0. $$

%
%
%
\beginsection 2.  Proof of the interpolation theorem

In this section we prove theorem 0.1.  The main construction is
contained in the following lemma.

\proclaim 2.1 Lemma:
Given $P_j$ and $m_j$ as in theorem 0.1, there
exists a discrete sequence $\{c_j\}$ contained in the $z_1$-axis and
$F \in \Aut\C^n$ such that for all $j$ we have
$$ F(z) = c_j + P_j(z-c_j) + O(|z-c_j|^{m_j+1}),\ \  z \to c_j. \eqno(2.1)$$
If in addition each polynomial $P_j$ satisfies (0.2), we may choose
$F \in \Aut_1\C^n$.   

\demo Proof:
{}For $j \ge 1$, let $\e_j = 2^{-j}$.  We will construct the sequence
$\{c_j\}$ inductively and construct $F$ as the limit of a
composition of automorphisms, each chosen inductively.

Let $K_0 = \emptyset$ and $K_1 = 2\overline{\B}$.  Let $c_1 = 3e_1$.  By
proposition 1.1 with $K = K_1$, there exists $\Psi_1 \in \Aut\C^n$
such that $\Psi_1(z) = c_1 + P_1(z-c_1) + O(|z-c_1|^{m_1+1})$, $z\to
c_1$ and $|\Psi_1(z)-z| + |\Psi_1^{-1}(z) - z| < \e_1$ for $z \in
K_1$.  If also $P_1$ satisfies (0.2), then we may choose $\Psi_1 \in 
\Aut_1\C^n$.  Let $F_0(z) = z$ and $F_1 = \Psi_1$.

\medskip
We will inductively choose automorphisms $\Psi_j$ and let $F_k =
\Psi_k \circ \cdots \circ \Psi_1$.  For the induction, suppose we have the
following.  

\ni
\item{(1)} Compact, convex sets $K_0 \subset K_1 \subset \cdots
\subset K_k$ with $j\B \cup F_{j-1}(j\B) \subset
K_j$ and 
$$\dist(K_{j-1} \cup F_{j-1}(j\B),\C^n\bs K_j) >\e_j$$
 for each $1 \le j \le k$.   
\item{(2)} Points $c_j \in K_{j+1}\bs K_j$ for $1 \le j \le k-1$ and
$c_k \in \C^n\bs K_k$ such that each $c_j$ is contained in the $z_1$-axis.
\item{(3)} $\Psi_j \in \Aut\C^n$ for $1 \le j \le k$ with
$$|\Psi_j(z)-z| + |\Psi_j^{-1}(z) - z| < \e_j$$
 for $z \in K_j$ (and $\Psi_j \in \Aut_1\C^n$ if $P_j$ satisfies (0.2)). 
\item{(4)} $F_k = \Psi_k \circ \Psi_{k-1} \circ \cdots \circ \Psi_1$
satisfying (2.1) for $1 \le j \le k$.  

\medskip
Given this, let $K_{k+1}$ be a compact, convex set in $\C^n$ with
$$(k+1)\B \cup K_k \cup F_k((k+1)\B) \cup \{c_k\} \subset K_{k+1}$$
and
$$\dist(K_k \cup F_k((k+1)\B),\C^n\bs K_{k+1}) >\e_{k+1}.$$
Since $K_{k+1}$ is compact, we can choose $c_{k+1}$ in the $z_1$-axis
so that $c_{k+1}, F_k(c_{k+1}) \in \C^n \bs K_{k+1}$.  Let $N =
\max\{m_1, \ldots, m_k\}+1$.   
Let $d_{k+1} = F_k(c_{k+1})$, and let $Q_{k+1} = P_{k+1}\circ
(F_k^{-1})_{m_{k+1},d_k}$.  Using proposition 1.1, choose 
$\Psi_{k+1} \in \Aut\C^n$ (with $\Psi_{k+1} \in \Aut_1\C^n$ if
$P_{k+1}$ satisfies (0.2) and $F_k \in \Aut_1\C^n$)) so that 

\ni
\item{(i)} $\Psi_{k+1}(z) = c_{k+1} + Q_{k+1}(z-d_{k+1}) + 
O(|z-d_{k+1}|^{m_{k+1}+1})$ as $z \to d_{k+1}$,
\item{(ii)} $\Psi_{k+1}(z) = z+O(|z-c_j|^N)$ as $z \to c_j$, $1\le j
\le k$,
\item{(iii)} $|\Psi_{k+1}(z) -z| + |\Psi_{k+1}^{-1}(z) - z| <\e_{k+1}$
for each $z \in K_{k+1}$. 

\medskip\ni
Taking $F_{k+1} = \Psi_{k+1} \circ F_k$, we obtain the induction
hypotheses at stage $k+1$.  

By proposition 1.2, the sequence $\{F_j\}$ converges uniformly on
compact subsets of $\Omega = \cup_{j=1}^\infty F_j^{-1}(K_j)$ to a
biholomorphic map $F$ from $\Omega$ onto $D = \cup_{j=1}^\infty K_j$.
Since $j\B \subset K_j$, we see that $D = \C^n$. 
Moreover, $F_j^{-1}(K_j) = F_{j-1}^{-1} \Psi_j^{-1}(K_j)$, and by
Rouch\'e's theorem [1, p.110] and induction hypotheses (1) and
(3), we see that $F_{j-1}(j\B) \subset \Psi_j^{-1}(K_j)$.  Hence $j\B
\subset F_j^{-1}(K_j)$, so also $\Omega =\C^n$.  

Hence $F \in \Aut\C^n$ satisfies (2.1) for all $j$ since each $F_k$
satisfies (2.1) for $1 \le j \le k$.  Finally, if each $P_j$ satisfies
(0.2), then each $\Psi_j \in \Aut_1\C^n$, so $F \in \Aut_1\C^n$.
\endpr

\demo Remark:
Given $R>0$, $1>\e>0$, we can replace $K_1$ by $(R+1)\overline{\B}$,
$c_1$ by $(R+2)e_1$, and $\e_j$ by $\e (\e_j/2)$ to construct a sequence
$\{c_j\}$ and $F$ satisfying the conclusions of the theorem and also
so that $|F(z) - z| + |F^{-1}(z) - z| < \e$ on $R\B$.   

\medskip
We next use a classical 1-variable interpolation result together with
theorem 0.1 to find an automorphism fixing each $je_1$ and having
prescribed jet there.  Our technique of using shears to map a given
discrete set in the $z_1$-axis to another was also used in the
proof of proposition 3.1 in [5].  However, for the current
application, we need the map to be tangent to a translation to
high order at each point in the discrete set. 

%
%
\proclaim 2.2 Corollary:
Let $P_j$ and $m_j$ be as in theorem 0.1. For each $R\ge 0$ and
$\e>0$ there exists $F \in \Aut\C^n$ such that
$$ \eqalign{ & F(z) = je_1 + P_j(z-je_1) + O(|z-je_1|^{m_j+1}),
                      \quad z \to je_1,\ j>R, \cr
   & |F(z)-z|+|F^{-1}(z)-z|<\e, \quad |z|<R. \cr}  $$
If each $P_j$ satisfies (0.2), then we can choose $F \in \Aut_1\C^n$.

\demo Proof: We shall first prove the corollary without
the last condition on $F$ (i.e., taking $R=0$).
By lemma 2.1, there exists $G \in \Aut\C^n$ and a sequence $\{c_j\}$
contained in the $z_1$-axis so that $G(z) = c_j + P_j(z-c_j) +
O(|z-c_j|^{m_j+1})$ as $z \to c_j$.  Also, $G \in \Aut_1\C^n$ if each
$P_j$ satisfies (0.2).  

To obtain the corollary we need only find $\Psi \in \Aut_1\C^n$ which
maps each $je_1$ to $c_j$ with $\Psi_{m_j, je_1}(z) = z-je_1$,
since then lemma 1.3 implies that $\Psi^{-1}\circ G\circ \Psi$
has the desired properties. 

To do this, let $\x_j \in \C$ such that $c_j = \x_j e_1$.
By a standard 1-variable interpolation result
[4, cor. 1.5.4], there exists an entire
function $f_1$ of one variable with $f_1(\z) = j + O(|\z-j|^{m_j+1})$
as $\z \to j$ for all $j \ge 1$.  Let $\Psi_1(z) = z+f_1(z_1)e_2$,
where $e_2=(0,1,0, \ldots, 0)$.  Then $\Psi_1 \in
\Aut_1\C^n$ and $\Psi_1(z) = j(e_1+e_2) + (z-je_1) +
O(|z-je_1|^{m_j+1})$ as $z \to je_1$. Thus $\Psi_1$ maps $je_1$ to
$j(e_1+e_2)$ and agrees with a translation to order $m_j+1$ at $je_1$.

Likewise, choosing $f_2$ entire with $f_2(\z) = \x_j - j +
O(|\z-j|^{m_j+1})$ as $\z \to j$, and taking $\Psi_2(z) =
z+f_2(z_2)e_1$, we see that $\Psi_2$ maps
$j(e_1+e_2)$ to $\x_j e_1 + j e_2$ and agrees with a translation to
order $m_j+1$ at $j(e_1+e_2)$.  Similarly, we can find $\Psi_3(z) = 
z_1+f_3(z_1)e_2$ so that $\Psi_3$ maps $\x_j e_1 + j e_2$ to $\x_j e_1
= c_j$ and agrees with a translation to order $m_j+1$ at
$\x_je_1+je_2$.

Let $\Psi = \Psi_3 \circ \Psi_2 \circ \Psi_1$.  Then $\Psi(z) = c_j +
(z-je_1) + O(|z-je_1|^{m_j+1})$ as $z \to je_1$, and by lemma 1.3,
we see that $\Psi^{-1}(z) = je_1 + (z-c_j) + O(|z-c_j|^{m_j+1})$ as
$z \to c_j$.  Also, $\Psi \in \Aut_1\C^n$ since each $\Psi_l$ has
Jacobian one. 

Hence taking $F = \Psi^{-1} \circ G \circ \Psi$, and applying lemma
1.3, we obtain the desired automorphism.

Our proof shows that the sequence $\{je_1\}$ could clearly be replaced
by any discrete sequence $\{d_j\}$ without repetition in the $z_1$-axis.
Moreover, to get the last condition on $F$, suppose that $\{d_j\}$
and $\{c_j\}$ lie outside the ball $R\overline{\B}$ (in the $z_1$ axis).
A simple argument using Runge's theorem and the Weierstrass theorem shows
that we can choose each $f_l$ as above and so that $|f_l|$ is small on
a neighborhood of the closed disk of radius $R$ in $\C$.  Hence by the
remark after the proof of lemma 2.1, we see that we may choose $F$
satisfying the conclusions of corollary 2.2 with $d_j$ in place of
$je_1$ and so that $|F(z) - z| + |F^{-1}(z) -z| < \e$ on $R\B$.
\endpr

\demo Proof of theorem 0.1:
Since $\{a_j\}$ and $\{b_j\}$ are tame, there exist $H_1, H_2 \in \Aut\C^n$
such that $H_1(a_j) = H_2(b_j) = j e_1$, and if the sequences are very
tame, then $H_1, H_2 \in \Aut_1\C^n$.  

Hence it suffices to construct $G \in \Aut\C^n$ (or $G\in \Aut_1\C^n$)
such that $G(je_1) = je_1$ and
$(H_2^{-1} G H_1)_{m_j,a_j} = P_j$.  By lemma 1.3, this latter
condition can be satisfied by making 
$$G_{m_j,je_1}(z) = (H_2)_{m_j, b_j} \circ P_j
\circ (H_1^{-1})_{m_j, je_1}(z) + O(|z|^{m_j+1}),\ z \to 0$$
for each $j$.  By corollary 2.2, we can find 
$G \in \Aut\C^n$ satisfying this condition and $G(je_1) = je_1$.  
Moreover, if each $P_j$ satisfies (0.2) and $H_1,
H_2 \in \Aut_1\C^n$, then we can choose $G \in \Aut_1\C^n$. 

{}Finally, taking $F = H_2^{-1} \circ G \circ H_1$, we obtain the
desired automorphism.  
\endpr

%
%
%
%
\beginsection 3.  Interpolation with approximation

We begin with the following lemma.

\proclaim 3.1 Lemma:
Let $H \in \Aut\C^n$, let $K$ be a compact, \pc\ set in $\C^n$,
and let $\{a_j\}_{j=1}^\infty$ be a discrete sequence disjoint
from $K$ and contained in the $z_1$-axis. Let $R>0$ with
$H(K) \subset R\B$ and let $\e>0$.  Then there exists
$\Psi \in \Aut\C^n$ such that $|\Psi(z) - H(z)|<\e$ on $K$,
$|\Psi^{-1}(z) - H^{-1}(z)| <\e$ on $H(K)$, and $\Psi(a_j) = (R+j)e_1$
for all $j$.  If $H \in \Aut_1\C^n$, then we can choose
$\Psi \in \Aut_1\C^n$.

\demo Proof:
Let $r>0$ with $K \subset r\B$.  Let $K_0 \subset r\B$ be compact and
\pc\ with $K \subset {\rm Int} K_0$, $K_0 \cap \{a_j\} = \emptyset$,
and $H(K_0) \subset R\B$.  Let $a_{j_1}, \ldots, a_{j_m}$ be the
points in $\{a_j\} \cap r\overline{\B}$, and let $\d>0$.  
Using the fact that the union of $K_0$ with finitely many points is again
\pc, we can apply proposition 1.1 $m$ times to find $\Psi_1 \in
\Aut_1\C^n$ such that $|\Psi_1(z) - z| + |\Psi_1^{-1}(z) - z| < \d$
for $z \in K_0$ and $\Psi_1(a_{j_k}) = H^{-1}((R+j_k)e_1)$ for $k = 1,
\ldots, m$.   

Let $\pi_2(z) = z_2$.  For fixed $v \in \C^n\bs\{0\}$ and $j\in \N$
we consider the 1-variable function
$g_j(\z) = \pi_2 H \Psi_1(a_j + \z v)$.  Since the
kernel of $D_{a_j} (\pi_2 \circ H \circ \Psi_1)$ is an $(n-1)$-dimensional
subspace for each $j$, we may choose $v$ arbitrarily near $e_2$
such that $g_j$ is nonconstant for each $j$ and hence so
that the image of $g_j$ omits at most one point in $\C$.
In particular, we may choose such a $v$ so that there exists
a $\C$-linear form $\lambda$ with $\lambda(v)=0$, $\lambda(e_1)=1$,
and $\lambda(a_j) \notin \lambda(r\overline{\B})$ if $j \notin
\{j_1, \ldots, j_m\}$. 

Choose $f$ entire on $\C$ such that $|f|<\d/2$ on $\lambda(r\B)$,
$f(\lambda(a_{j_k})) = 0$ for each $k=1, \ldots, m$, and $|\pi_2 H
\Psi_1(a_j + f(\lambda(a_j)) v)| = R+j$ for each $j \notin \{j_1,
\ldots, j_m\}$.   Let $\Psi_2(z) = z + f(\lambda(z)) v$.  Then $\Psi_2
\in \Aut_1\C^n$, 
$|\Psi_2(z)-z| + |\Psi_2^{-1}(z)-z| < \d$ on $r\B$, $H
\Psi_1 \Psi_2(a_{j_k}) = (R+j_k)e_1$ for $k=1, \ldots, m$, and $|\pi_2 H
\Psi_1 \Psi_2(a_j)| = R+j$ for $j \notin \{j_1, \ldots, j_m\}$. 

Using a composition of two shears as in the proof of corollary
2.2, we can find $\Psi_3 \in \Aut_1\C^n$ such that $|\Psi_3(z)-z| +
|\Psi_3^{-1}(z)-z| < \d$ on $R\B$, $\Psi_3((R+j_k)e_1) = (R+j_k)e_1$
for $k=1, \ldots, m$, and $\Psi_3 H \Psi_1 \Psi_2(a_j) = (R+j)e_1$ for
each $j \notin \{j_1, \ldots, j_m\}$. 

Let $\Psi = \Psi_3 \circ H \circ \Psi_1 \circ \Psi_2$.  Then
$\Psi(a_j) = (R+j)e_1$ for all $j$, and for $\d$ sufficiently small,
we have $|\Psi(z) - H(z)| < \e$ on $K$ and $|\Psi^{-1}(z) - H^{-1}(z)|
<\e$ on $H(K)$, so the lemma follows.  Since each $\Psi_l \in
\Aut_1\C^n$, we see that $\Psi \in \Aut_1\C^n$ if $H \in \Aut_1\C^n$.
\endpr

\demo Proof of theorem 0.2:
Since the sequences $\{a_j\}$ and $\{b_j\}$ are tame, there
exist automorphisms $H_1$, $H_2$ of $\C^n$ such that
$H_1(a_j)=je_1$ and $H_2(b_j)=je_1$ for each $j$. Replacing
$\Phi$ by $H_2\circ \Phi \circ H_1^{-1}$,  $K$ by
$H_1(K)$, and adjusting the jets $P_j$ as in the proof of theorem 0.1,
we reduce the problem to the case when 
$a_j=b_j=je_1$ for all $j$.

Choose a larger \pc\ set $L\subset \C^n\bs \{je_1 \colon j\in \N\}$
such that $K\subset {\rm Int}\, L$.  Fix $\eta>0$, and 
choose an integer $r>0$ such that
$L\cup \Phi(L) \subset r\B$. By lemma 3.1 there exist
automorphisms $\Psi, \Theta\in\Aut\C^n$ such that
$$ \eqalign{ & \Psi(je_1)=\Theta(je_1)= (r+j)e_1,
               \quad j=1,2,3,\ldots, \cr
   & |\Psi(z)-\Phi(z)|<\eta, \quad z\in L, \cr
   & |\Theta(z)-z| + |\Theta^{-1}(z)-z| <\eta, \quad z\in \Phi(L).\cr} $$
By corollary 2.2, there exists $G \in \Aut\C^n$ such that
$G((r+j)e_1)=(r+j)e_1$ for $j \in \N$,
%
%
$|G(z)-z|<\eta$ for $|z|\le r$, and such that for each $j\in \N$
the jet $Q_j=G_{m_j,(r+j)e_1}$ satisfies
$$ P_j(z)= (\Theta^{-1})_{m_j,(r+j)e_1} \circ Q_j \circ
           \Psi_{m_j,je_1}(z) + O(|z|^{m_j+1}),\quad z\to 0. $$

Let $F=\Theta^{-1}\circ G\circ \Psi$.
Then $F(je_1)=je_1$ for each $j\ge 1$, $F$ satisfies (0.3)
provided that $\eta>0$ is chosen sufficiently small
(depending on $\e$ and $\dist(K,\C^n\bs L)$), and
lemma 1.3 shows that $F$ satisfies (0.1), with
$a_j=b_j=je_1$. 

Finally, if the sequences $a_j$ and $b_j$ are very tame,
the $P_j$'s satisfy (0.2), and $\Phi \in\Aut_1\C^n$, then
the automorphisms $\Psi$, $\Theta$, and $G$ can be chosen
with Jacobian one.
\endpr

%
%
%
%

\bigskip
\ni\bf References. \rm
\medskip

\item{1.} Chirka, E.: \it Complex Analytic Sets. \rm
Kluwer, Dordrecht 1989

\item{2.} Forstneric, F.:
Holomorphic automorphism groups of $\C^n$: A survey.
(The Proceedings `Complex Analysis and Geometry', Ed.\
V.\ Ancona, E.\ Ballico, A.\ Silva; pp.\ 173--200)
\it Lecture Notes in Pure and Applied Mathematics \bf 173\rm,
Marcel Dekker, New York 1996

\item{3.} Forstneric, F.: Interpolation by holomorphic automorphisms
and embeddings in $C^n$.  Preprint, 1996

\item{4.} H\"ormander, L.:
An Introduction to Complex Analysis in Several Variables,
\rm 3rd ed. Amsterdam: North Holland 1990

\item{5.} Rosay, J.-P., Rudin, W.:
Holomorphic maps from $\C^n$ to $\C^n$.
Trans.\ Amer.\ Math.\ Soc.\ {\bf 310}, 47--86 (1988)

\bye